\documentclass[a4paper,twoside,10pt,reqno]{article}
\usepackage[english]{babel}
\usepackage[dvips]{graphicx}
\usepackage[T1]{fontenc}
\usepackage[latin1]{inputenc}
\usepackage{amssymb,amsthm,dsfont,mathrsfs}
\usepackage{amsfonts,amsmath,euscript}
\usepackage{color}

\newtheorem{defi}{Definition}[section]
\newtheorem{teo}{Theorem}[section]
\newtheorem{pro}[teo]{Proposition}
\newtheorem{cor}[teo]{Corollary}

\newtheorem{rem}{Remark}[section]
\newtheorem{ex}{Example}[section]

\newcounter{example}[section]

           \addto{\captionsenglish}{}

\newcommand{\hs}{\hspace{3pt}}
\newcommand{\dst}{\displaystyle}
\newcommand{\dem}{{\bf Dem. }}
\newcommand{\fdem}{$\square$}

\newcommand{\nn}{\nonumber}

\newcommand{\nid}{\noindent}

\newcommand{\titulo}[1]{\mbox{} \\ \noindent \textit{\textbf{\Large #1}}\\}

\renewcommand{\abstract}[1]{{\small \noindent \textbf{Abstract:} #1\\}}

\newcommand{\pchave}[1]{{\small \noindent \textbf{Keywords:} #1\\}}

\begin{document}

\begin{center}
\titulo{Extremal dependence: some contributions}
\end{center}

\vspace{0.5cm}

\textbf{Helena Ferreira} Department of Mathematics, University of
Beira
Interior, Covilhã, Portugal\\

\textbf{Marta Ferreira} Department of Mathematics, University of
Minho, Braga, Portugal\\

\abstract{Due to globalization and relaxed market regulation, we
have assisted to an increasing of extremal dependence in
international markets. As a consequence, several measures of tail
dependence have been stated in literature in recent years, based on
multivariate extreme-value theory. In this paper we present a tail
dependence function and an extremal coefficient of dependence
between two random vectors that extend existing ones. We shall see
that in weakening the usual required dependence allows to assess the
amount of dependence in $d$-variate random vectors based on
bidimensional techniques. Very simple estimators will be stated and
can be applied to the well-known \emph{stable tail dependence
function}. Asymptotic normality and strong consistency will be
derived too. An application to financial markets will be presented
at the end.}

\pchave{multivariate extreme value theory, tail dependence, extremal
coefficients}

\section{Introduction}
Dependence between extremal events have increased in recent time
periods in financial markets, especially during bear markets and
market crashes. The globalization and the lack of supervision are
well-known contributions for this phenomena. Therefore, modern risk
management is highly interested in assessing the amount of extremal
dependence.
The concept of tail dependence is
the current tool used to this end, although it was first introduced
far back in the sixties (Sibuya \cite{sib}, 1960; Tiago de Oliveira
\cite{tiago}, 1962/63). Tail dependence coefficients measure the
probability of occurring extreme values
 for one random variable (r.v.) given that
another assumes an extreme value too. These coefficients can be
defined via copulas of random vectors which refers to their
dependence structure concerning extreme events independently of
their marginal distributions. The tail dependence coefficient,
\begin{eqnarray}\label{tdc}
\dst\lambda=\lim_{t\downarrow 0}P(F_X(X)>1-t|F_Y(Y)>1-t),
\end{eqnarray}
where $F_X$ and $F_Y$ are the distribution functions (d.f.'s) of $X$
and $Y$, respectively, is perhaps the most referred in literature
and characterizes the dependence in the tail of a random pair
$(X,Y)$, i.e., $\lambda>0$ corresponds to tail dependence and
$\lambda=0$ means tail independence. There are several references on
this topic (besides the two above) and thus we point out only some
of them: Ledford and Tawn (\cite{led+tawn1,led+tawn2}, 1996, 1997),
Joe (\cite{joe}, 1997), Coles \emph{et al.} (\cite{coles+}, 1999),
Embrechts \emph{et al.} (\cite{emb+03}, 2003).

Multivariate formulations for tail dependence coefficients can be
used to describe the amount of dependence in the orthant tail of a
multivariate distribution (Marshall-Olkin \cite{mar+olkin} 1967;
Wolff \cite{wolff} 1980; Nelsen \cite{nelsen} 1996; Schmid and
Schmidt \cite{schmid+schmidt} 2007; Li \cite{li1,li2,li3} 2006,
2008, 2009, among others). These have been increasingly used in the
most recent and higher demanding times. Most of the multivariate
measures consider that extremal events must occur to all the
components of the random vector, and obviously they are more
complicated to deal with and to understand than in the bivariate
case. Not surprisingly, applications hardly go any further than the
three-dimensional case.

But maybe this is a too demanding condition and the occurrence of at
least one extremal event in sub-vectors of a random vector can be
enough to assess dependence. As already mentioned, financial markets
are increasingly connected and the occurrence of at least one market
crash, for instance, in Europe, will certainly influence a negative
behavior in USA markets.

Based on this, we define a new tail dependence function for a random
vector as a measure of the probability of occurring extreme values
 for the maximum of one sub-vector  given that
the maximum of another assumes an extreme value too. At the unit
point, this function gives rise to the here called \emph{extremal
coefficient of dependence} as it relates to the well-known extremal
coefficient (Tiago de Oliveira 1962-63, Smith 1990). These extend,
respectively, the concept of \emph{upper tail dependence function}
and \emph{upper tail dependence coefficient} already stated in
literature (see Schmidt and Stadtm\"{u}ller (2006) and references
therein). In deriving the moments of the random variables involved
in this approach, we find very simple estimators that can be also
applied to the well-known \emph{stable tail dependence function}
(for a survey on this function see e.g. Beirlant \emph{et al.}
\cite{beirl+}, 2004). Asymptotic normality and strong consistency
are proved.\\

This paper is organized as follows. In Section \ref{sed} we define
our new \emph{upper-tail dependence function} and the \emph{extremal
coefficient of dependence}. We present some properties and examples.
We also analyze the case of asymptotic independence. In Section
\ref{sestim} we present estimators and derive the respective
properties of asymptotic normality and strong consistency. Section
\ref{saplic} illustrates our approach through an application to
financial data.

\section{Extremal dependence between two random vectors}\label{sed}

Let $\mathbf{X}=(X_1,...,X_d)$ be a random vector with d.f. $F$ and
continuous marginal d.f.'s $F_i$. For $I\subset\{1,...,d\}$, define
$M(I)=\bigvee_{i\in I}F_i(X_i)$ and $\mathbf{X}_I$ the sub-vector of
$\mathbf{X}$ having r.v.'s with indexes in $I$. Consider $C_F$ the
copula function of $F$, i.e.,
\begin{eqnarray}\label{copula}
F(x_1,...,x_d)=C_F(F_1(x_1),...,F_d(x_d)), \,\,\,(x_1,...,x_d)\in \mathbb{R}^d.
\end{eqnarray}

We are going to study the dependence between extremal events
concerning two sub-vectors, $\mathbf{X}_{I_1}$ and
$\mathbf{X}_{I_2}$, where $I_1$ and $I_2$ are disjoint subsets of
$\{1,...,d\}$.

We start by extending in Definition \ref{dutdf} the concept of
\emph{upper tail dependence function} (see Schmidt and
Stadtm\"{u}ller (2006) and references therein) and from this we
define a new tail dependence coefficient between two random vectors.

\begin{defi}\label{dutdf}
Let $I_1$ and $I_2$ be two non-empty subsets of $\{1,...,d\}$. The
upper-tail dependence function of $\mathbf{X}_{I_1}$ given
$\mathbf{X}_{I_2}$ is defined as, for $(x,y)\in(0,\infty)^2$,
\begin{eqnarray}\label{utdfcond}
\Lambda_{U}^{(I_1|I_2)}(x,y)=\lim_{t\to\infty}P\Big(M(I_1)>1-\frac{x}{t} \Big| M(I_2)>1-\frac{y}{t}\Big),
\end{eqnarray}
provided the limit exists.
\end{defi}

By taking $x=y=1$, we have
\begin{eqnarray}\label{utdcoefcond}
\Lambda_{U}^{(I_1|I_2)}(1,1)=\lim_{t\to\infty}P\Big(M(I_1)>1-\frac{1}{t} \Big| M(I_2)>1-\frac{1}{t}\Big),
\end{eqnarray}
which is a tail dependence coefficient greater than the one
considered in Li and Sun (2008),
\begin{eqnarray}\label{lisun}
\gamma=\lim_{t\to\infty}P\Big(\bigcap_{i\in {I_1}}F_i(X_i)>1-\frac{1}{t} \Big| \bigcup_{i\in {I_2}}F_i(X_i)>1-\frac{1}{t}\Big),
\end{eqnarray}
which in turn is greater than the coefficient of  Li (2009) for
$I_1=\{1,...,d\}-I_2$,
\begin{eqnarray}\label{schmistad}
\tau=\lim_{t\to\infty}P\Big(\bigcap_{i\in {I_1}}F_i(X_i)>1-\frac{1}{t} \Big| \bigcap_{i\in {I_2}}F_i(X_i)>1-\frac{1}{t}\Big).
\end{eqnarray}

The tail dependence coefficient $\Lambda_{U}^{(I_1|I_2)}(1,1)$ give
us information about the probability of occurring some extreme value
in $\{F_i(X_i), i\in I_1\}$ given that some extreme value occurs in
$\{F_i(X_i), i\in I_2\}$.\\

Before presenting the properties of function
$\Lambda_{U}^{(I_1|I_2)}(x,y)$ that will be the basis for the
definition of our coefficient, consider the following notation:

for $(x,y)\in (0,\infty)^2$, $\emptyset\subseteq
I_1,I_2\subseteq\{1,...,d\}$ and $i\in \{1,...,d\}$, let
\begin{eqnarray}\label{}
a_i^{(I_1,I_2)}(x,y)=x\mathbf{1}_{I_1}(i)+y\mathbf{1}_{I_2}(i)+\infty\mathbf{1}_{\overline{I_1\cup I_2}}(i),
\end{eqnarray}
where $\mathbf{1}(\cdot)$ is the indicator function, and
\begin{eqnarray}\label{}
l^{(I_1,I_2)}(x,y)=-\log F(a_1^{(I_1,I_2)}(x,y),...,a_d^{(I_1,I_2)}(x,y)),
\end{eqnarray}
with the convention that, when some of the arguments of $F$ are
$\infty$ we understand the limit of $F$ as those arguments tend to
$\infty$.

If $F$ is a multivariate extreme value distribution (MEV) with unit
Fr\'{e}chet marginals, we have
\begin{eqnarray}\label{}
l^{(I_1,I_2)}(x^{-1},x^{-1})=-\log (\exp(-x))^{\epsilon_{I_1\cup I_2}}=x\epsilon_{I_1\cup I_2},
\end{eqnarray}
where $\epsilon_{I_1\cup I_2}$ is the extremal coefficient of
$\mathbf{X}_{I_1\cup I_2}$ (Tiago de Oliveira 1962-63, Smith 1990).

\begin{pro}\label{p.1}
If $F$ is a MEV distribution with unit Fr\'{e}chet marginals, then
function $\Lambda_{U}^{(I_1|I_2)}(x,y)$ is defined and verifies
\begin{eqnarray}\label{}
\Lambda_{U}^{(I_1|I_2)}(x,y)=1+\frac{x\epsilon_{I_1}}{y\epsilon_{I_2}}-\frac{l^{(I_1,I_2)}(x^{-1},y^{-1})}{y\epsilon_{I_2}},
\end{eqnarray}
\end{pro}
\dem We have
\begin{eqnarray}\label{p.1.1}
\Lambda_{U}^{(I_1|I_2)}(x,y)=\lim_{t\to\infty}1+\frac{1-P\big(M(I_1)\leq 1-\frac{x}{t}\big)}{1-P\big(M(I_2)\leq 1-\frac{y}{t}\big)}
-\frac{1-P\big(M(I_1)\leq 1-\frac{x}{t},M(I_2)\leq 1-\frac{y}{t}\big)}{1-P\big(M(I_2)\leq 1-\frac{y}{t}\big)}.
\end{eqnarray}
On the other hand, for
$\alpha_i^{(I_1,I_2)}(u,v)=u\mathbf{1}_{I_1}(i)+v\mathbf{1}_{I_2}(i)+\mathbf{1}_{\overline{I_1\cup
I_2}}(i)$, it holds
\begin{eqnarray}\label{relhuang}
\begin{array}{rl}
&{\dst\lim_{t\to\infty}}-t\log P\big(M(I_1)\leq 1-\frac{x}{t},M(I_2)\leq 1-\frac{y}{t}\big)
\vspace{0.35cm}\\
=&{\dst\lim_{t\to\infty}}-t\log C_F\Big(\alpha_1^{(I_1,I_2)}\big(1-\frac{x}{t},1-\frac{y}{t}\big),...,
\alpha_d^{(I_1,I_2)}\big(1-\frac{x}{t},1-\frac{y}{t}\big)\Big)\vspace{0.35cm} \\
=&{\dst\lim_{t\to\infty}}-\log C_F\Big(\alpha_1^{(I_1,I_2)}\big(\big(1-\frac{x}{t}\big)^t,\big(1-\frac{y}{t}\big)^t\big),...,
  \alpha_d^{(I_1,I_2)}\big(\big(1-\frac{x}{t}\big)^t,\big(1-\frac{y}{t}\big)^t\big)\Big)\vspace{0.35cm} \\
=&-\log C_F\Big(\alpha_1^{(I_1,I_2)}\big(\exp(-x),\exp(-y)\big),...,
  \alpha_d^{(I_1,I_2)}\big(\exp(-x),\exp(-y)\big)\Big)\vspace{0.35cm} \\
=&-\log F\Big(  a_1^{(I_1,I_2)}\big(x^{-1},y^{-1}\big),...,a_d^{(I_1,I_2)}\big(x^{-1},y^{-1}\big)\Big)\vspace{0.35cm} \\
=&l^{(I_1,I_2)}\big(x^{-1},y^{-1}\big).
\end{array}
\end{eqnarray}
Therefore, dividing the numerator and denominator of the fractions
in (\ref{p.1.1}) by $t$, we obtain
\begin{eqnarray}\label{}
\begin{array}{rl}
\Lambda_{U}^{(I_1|I_2)}(x,y)=&\dst1+\frac{l^{(I_1,\emptyset)}(x^{-1},x^{-1})}{l^{(\emptyset,I_2)}(y^{-1},y^{-1})}
-\frac{l^{(I_1,I_2)}(x^{-1},y^{-1})}{l^{(\emptyset,I_2)}(y^{-1},y^{-1})}\vspace{0.35cm} \\
=&\dst1+\frac{-\log (\exp(-x))^{\epsilon_{I_1}}}{-\log (\exp(-y))^{\epsilon_{I_2}}}
-\frac{l^{(I_1,I_2)}(x^{-1},y^{-1})}{-\log (\exp(-y))^{\epsilon_{I_2}}}. \,\,$\fdem$\\
\end{array}
\end{eqnarray}

\vspace{0.5cm}

Therefore, under the conditions of Proposition \ref{p.1}, we have
$$
y{\epsilon_{I_2}}\Lambda_{U}^{(I_1|I_2)}(x,y)=x{\epsilon_{I_1}}\Lambda_{U}^{(I_2|I_1)}(y,x)
=x{\epsilon_{I_1}}+y{\epsilon_{I_2}}-l^{(I_1,I_2)}(x^{-1},y^{-1})
$$
and we will denote this common value as
$\Lambda_{U}^{(I_1,I_2)}(x,y)$.
\begin{defi}\label{dutdf}
The upper-tail dependence function for random vector
$(\mathbf{X}_{I_1},\mathbf{X}_{I_2})$ with d.f. MEV and unit
Fr\'{e}chet marginals is defined as
\begin{eqnarray}\label{utdf}
\Lambda_{U}^{(I_1,I_2)}(x,y)=x\epsilon_{I_1}+y\epsilon_{I_2}-l^{(I_1,I_2)}(x^{-1},y^{-1})
\end{eqnarray}
and the extremal coefficient of dependence between
$\mathbf{X}_{I_1}$ and $\mathbf{X}_{I_2}$ is given by
$\Lambda_{U}^{(I_1,I_2)}(1,1)$, which we denote
$\epsilon_{(I_1,I_2)}$ and hence
\begin{eqnarray}\label{extcoef}
\epsilon_{(I_1,I_2)}=\epsilon_{I_1}+\epsilon_{I_2}-\epsilon_{I_1\cup I_2}.
\end{eqnarray}
\end{defi}

The upper-tail dependence function (\ref{utdf}) generalizes the
relation of Huang (1992) corresponding to $I_1=\{1\}$ and
$I_2=\{2\}$,
\begin{eqnarray}\label{lambdaf}
\Lambda_{U}(x,y)=x+y-l^H_{(F_1(X_1),F_2(X_2))}(x,y),
\end{eqnarray}
where the stable tail dependence function in the right-side is given
by
\begin{eqnarray}\label{stf}
l^H_{(F_1(X_1),F_2(X_2))}(x,y)
=\lim_{t\to\infty}tP\Big(F_1(X_1)>1-\frac{x}{t} \vee F_2(X_2)>1-\frac{y}{t}\Big).
\end{eqnarray}
Observe that by (\ref{p.1.1}) we also obtain
\begin{eqnarray}\label{stdfnossa}
l^{(I_1,I_2)}(x^{-1},y^{-1})=\lim_{t\to\infty}tP\Big(M(I_1)>1-\frac{x}{t} \vee M(I_2)>1-\frac{y}{t}\Big)=l^H_{(M(I_1),M(I_2))}(x,y).
\end{eqnarray}
Moreover, the upper-tail dependence function in (\ref{utdf}) can be
can be viewed as an extension of the bivariate upper-tail dependence
function of Schmidt and Stadtm\"{u}ller (\cite{schmidt+stadt},
2006), defined as
\begin{eqnarray}\label{utdfS}
\Lambda^S_{(F_1(X_1),F_2(X_2))}(x,y)=\lim_{t\to\infty}t
P\Big(F_1(X_1)>1-\frac{x}{t},F_2(X_2)>1-\frac{y}{t}\Big),
\end{eqnarray}
 by taking in this limit the random pair $(M(I_1),M(I_2))$ instead of $(F(X_1),F(X_2))$. 
At the unit vector, the Schmidt and Stadtm\"{u}ller upper-tail
dependence function corresponds to the tail dependence coefficient
$\lambda$ in (\ref{tdc}), i.e.,
$\lambda=\Lambda^S_{(F_1(X_1),F_2(X_2))}(1,1)$.

In the following, we present the expression of the tail-dependence
function $\Lambda_U^{(I_1,I_2)}(x,y)$ and the value of the
corresponding extremal coefficient $\epsilon_{(I_1,I_2)}$ for a
$d$-variate random vector $\mathbf{X}$ with well-known distribution
functions for its margins.\\

\begin{ex}
Consider vector $\mathbf{X}$ with unit Fr\'{e}chet margins and
copula function
$C_{\mathbf{X}}(u_1,...,u_d)=\prod_{l=1}^\infty\prod_{k=-\infty}^\infty
u_1^{\alpha_{lk1}}\wedge ...\wedge u_d^{\alpha_{lkd}}$, where
$u_j\in [0,1]$, $j=1,...,d$, and $\{\alpha_{lkj}, -\infty
<k<\infty,1\leq j\leq d,l\geq 1\}$ is a family of non negative
constants such that $\sum_{l=1}^\infty\sum_{k=-\infty}^\infty
\alpha_{lkj}=1$, $j=1,...,d$. The distribution of $\mathbf{X}$ is
the MEV marginal distribution of \emph{multivariate maxima of moving
maxima} processes considered in Smith and Weissman
(\cite{smith+weissman}, 1996). We have
$$
l^{(I_1,I_2)}(x,y)=-\log C\Big(e^{-a_1^{-1}(x,y)},...,e^{-a_d^{-1}(x,y)}\Big)
=\sum_{l=1}^\infty\sum_{k=-\infty}^\infty\bigvee_{j=1}^d a_j^{-1}(x,y)\alpha_{lkj}.
$$
Therefore,
$$
\begin{array}{rl}
&\dst\Lambda_U^{(I_1,I_2)}(x,y)=l^{(I_1,\emptyset)}(x^{-1},x^{-1})+l^{(\emptyset,I_2)}(y^{-1},y^{-1})-l^{(I_1,I_2)}(x^{-1},y^{-1}) \vspace{0.35cm}\\
=&\dst x \sum_{l=1}^\infty\sum_{k=-\infty}^\infty\bigvee_{j\in I_1}\alpha_{lkj}
+y\sum_{l=1}^\infty\sum_{k=-\infty}^\infty\bigvee_{j\in I_2}\alpha_{lkj}
-\sum_{l=1}^\infty\sum_{k=-\infty}^\infty\bigg(\Big(x\bigvee_{j\in I_1}\alpha_{lkj} \Big)\vee\Big(y\bigvee_{j\in I_2}\alpha_{lkj} \Big) \bigg)
\end{array}
$$
and
$$
\begin{array}{c}
\dst\epsilon_{(I_1,I_2)}(x,y)=\dst  \sum_{l=1}^\infty\sum_{k=-\infty}^\infty\bigvee_{j\in I_1}\alpha_{lkj}
+\sum_{l=1}^\infty\sum_{k=-\infty}^\infty\bigvee_{j\in I_2}\alpha_{lkj}
-\sum_{l=1}^\infty\sum_{k=-\infty}^\infty\bigvee_{j\in I_1\cup I_2}\alpha_{lkj}
\end{array}
$$
Illustrating with
$$
C_{\mathbf{X}}(u_1,u_2,u_3,u_4)=(u_1^{1/8}\wedge u_2^{1/8}\wedge u_3^{1/8}\wedge
u_4^{1/8}).(u_1^{5/8}\wedge u_2^{4/8}\wedge u_3^{7/8}\wedge
u_4^{1/8}).(u_1^{1/8}\wedge u_2^{2/8}).(u_1^{1/8}\wedge
u_2^{1/8}\wedge u_4^{6/8}),
$$
$I_1=\{1,2\}$ and $I_2=\{3,4\}$, we obtain
$$
\begin{array}{rl}
\Lambda_U^{(I_1,I_2)}(x,y)=& \Big( \frac{1}{8}+\frac{5}{8}+\frac{2}{8}+\frac{1}{8}\Big)x
+ \Big( \frac{1}{8}+\frac{7}{8}+\frac{6}{8}\Big)y   
-\Big( \big(x\frac{1}{8}\vee y\frac{1}{8}\big)+\big(x\frac{5}{8}\vee y\frac{7}{8}\big)
+\big(x\frac{2}{8}\big)+\big(x\frac{1}{8}\vee y\frac{6}{8}\big)\Big)\vspace{0.35cm}\\
=& \frac{9}{8}x+\frac{14}{8}y- \Big( \big(x\frac{1}{8}\vee y\frac{1}{8}\big)+\big(x\frac{5}{8}\vee y\frac{7}{8}\big)
+x\frac{2}{8}+\big(x\frac{1}{8}\vee y\frac{6}{8}\big)\Big)
\end{array}
$$
and
$$
\begin{array}{c}
\epsilon_{(I_1,I_2)}=\frac{9}{8}+\frac{14}{8}- \big(\frac{1}{8}+ \frac{7}{8}+\frac{2}{8}+\frac{6}{8}\big)
=\frac{7}{8}
\end{array}
$$
Similarly, if $I_1=\{1,2\}$ and $I_2=\{4\}$ we obtain
$$
\begin{array}{rl}
\Lambda_U^{(I_1,I_2)}(x,y)=\frac{9}{8}x+y- \Big( \big(x\frac{1}{8}\vee y\frac{1}{8}\big)+\big(x\frac{5}{8}\vee y\frac{1}{8}\big)
+x\frac{2}{8}+\big(x\frac{1}{8}\vee y\frac{6}{8}\big)\Big)
\end{array}
$$
and
$$
\begin{array}{c}
\epsilon_{(I_1,I_2)}=\frac{9}{8}+1- \big(\frac{1}{8}+ \frac{5}{8}+\frac{2}{8}+\frac{6}{8}\big)
=\frac{3}{8} .
\end{array}
$$
\end{ex}

\begin{ex}
For the symmetric logistic model we have
$$
l^{(I_1,I_2)}(x,y)=-\log F(a_1^{(I_1,I_2)}(x,y),...,a_d^{(I_1,I_2)}(x,y))=\Big(\sum_{j=1}^d (a_j^{(I_1,I_2)}(x,y))^{-1/\theta}\Big)^{\theta}
$$
with $\theta\in(0,1]$, $x,y>0$. Therefore,
$$
\begin{array}{rl}
\dst\Lambda_U^{(I_1,I_2)}(x,y)=&l^{(I_1,\emptyset)}(x^{-1},x^{-1})+l^{(\emptyset,I_2)}(y^{-1},y^{-1})-l^{(I_1,I_2)}(x^{-1},y^{-1}) \vspace{0.35cm}\\
=&\dst \Big(\sum_{j\in I_1}x^{1/\theta}\Big)^{\theta} +\Big(\sum_{j\in I_2}y^{1/\theta}\Big)^{\theta}
-\Big(\sum_{j\in I_1}x^{1/\theta}+\sum_{j\in I_2}y^{1/\theta}\Big)^{\theta}     \vspace{0.35cm}\\
=& |I_1|^{\theta}x+|I_2|^{\theta}y-\big(|I_1|x^{1/\theta}+|I_2|y^{1/\theta}\big)^{\theta}
\end{array}
$$
and
$$
\begin{array}{c}
\epsilon_{(I_1,I_2)}=|I_1|^{\theta}+|I_2|^{\theta}-\big(|I_1|+|I_2|\big)^{\theta}.\\
\end{array}
$$
\end{ex}
\vspace{0.5cm}

\begin{pro}\label{pbounds}
Under the conditions of Proposition \ref{p.1} 
we have
\begin{enumerate}
\item[(i)]\label{enq1} $0\leq\Lambda_U^{(I_1,I_2)}(x,y)\leq x\epsilon_{I_1}\wedge
y\epsilon_{I_2}$\\
\item[(ii)]\label{enq2} $0\leq\epsilon_{(I_1,I_2)}\leq \epsilon_{I_1}\wedge
\epsilon_{I_2}$.
\end{enumerate}
\end{pro}
\dem
\begin{enumerate}
\item[(i)]
 The left inequality is straightforward by the definition of
$\Lambda_U^{(I_1,I_2)}(x,y)$ in (\ref{utdf}). Observe also that,
since $\mathbf{X}$ has MEV distribution, it is associated (in
the sense of Joe \cite{joe}, 1997; Theorem 6.7) and hence, for
all $I_1,I_2\subset\{1,...,d\}$, 
$$
\begin{array}{rl}
&F(a_1^{(I_1,I_2)}(x^{-1},y^{-1}),...,a_d^{(I_1,I_2)}(x^{-1},y^{-1}))\vspace{0.35cm}\\
\geq &
F(a_1^{(I_1,\emptyset)}(x^{-1},x^{-1}),...,a_d^{(I_1,\emptyset)}(x^{-1},x^{-1}))
F(a_1^{(\emptyset,I_2)}(y^{-1},y^{-1}),...,a_d^{(\emptyset,I_2)}(y^{-1},y^{-1})),
\end{array}
$$
leading to the same conclusion, i.e.,
$$
l^{(I_1,\emptyset)}(x^{-1},x^{-1})+l^{(\emptyset,I_2)}(y^{-1},y^{-1})-l^{(I_1,I_2)}(x^{-1},y^{-1})\geq 0.
$$
On the other hand,
$$
\begin{array}{rl}
&F(a_1^{(I_1,I_2)}(x^{-1},y^{-1}),...,a_d^{(I_1,I_2)}(x^{-1},y^{-1}))\vspace{0.35cm}\\
\leq &
F(a_1^{(I_1,\emptyset)}(x^{-1},x^{-1}),...,a_d^{(I_1,\emptyset)}(x^{-1},x^{-1}))\wedge
F(a_1^{(\emptyset,I_2)}(y^{-1},y^{-1}),...,a_d^{(\emptyset,I_2)}(y^{-1},y^{-1}))
\end{array}
$$
and hence
$$
\begin{array}{rl}
&l^{(I_1,\emptyset)}(x^{-1},x^{-1})+l^{(\emptyset,I_2)}(y^{-1},y^{-1})-l^{(I_1,I_2)}(x^{-1},y^{-1})\vspace{0.35cm}\\
\leq &
l^{(I_1,\emptyset)}(x^{-1},x^{-1})+l^{(\emptyset,I_2)}(y^{-1},y^{-1})-
\big(l^{(I_1,\emptyset)}(x^{-1},x^{-1}) \vee l^{(\emptyset,I_2)}(y^{-1},y^{-1})\big)   \vspace{0.35cm}\\
= &    x\epsilon_{I_1}\wedge
y\epsilon_{I_2}.\,\,  $\fdem$
\end{array}
$$
\end{enumerate}
The result in (i) agrees with the one for the bivariate case.
Observe that, from the proof above we can also conclude that the
boundary cases correspond to,
respectively, independence and total dependence.\\

\begin{rem}\label{rground}
With the conventions $1/0:=\infty$ and $1/\infty:=0$, we can define
$\Lambda_U^{(I_1,I_2)}(x,y)$ in
$[0,\infty]^2\backslash\{(\infty,\infty)\}$ and found
$\Lambda_U^{(I_1,I_2)}(0,y)=0=\Lambda_U^{(I_1,I_2)}(x,0)$,
$\Lambda_U^{(I_1,I_2)}(\infty,y)=y\epsilon_{I_2}$ and
$\Lambda_U^{(I_1,I_2)}(x,\infty)=x\epsilon_{I_1}$.
\end{rem}

\begin{pro}\label{pderiv}
Under the conditions of Proposition \ref{p.1} and Remark
\ref{rground}, for each $y\geq 0$, the partial derivative
$\partial\Lambda_U^{(I_1,I_2)}/
\partial x$ exists for almost all $x>0$, and 
$$
\begin{array}{c}
0\leq \frac{\partial}{\partial x}\Lambda_U^{(I_1,I_2)}(x,y)\leq |I_1|.
\end{array}
$$
Similarly, for each $x\geq 0$, the partial derivative
$\partial\Lambda_U^{(I_1,I_2)}/
\partial x$ exists for almost all $y>0$, and 
$$
\begin{array}{c}
0\leq \frac{\partial}{\partial y}\Lambda_U^{(I_1,I_2)}(x,y)\leq |I_2|.
\end{array}
$$
Also, the functions $x\mapsto\partial\Lambda_U^{(I_1,I_2)}(x,y)/
\partial y$ and $y\mapsto\partial\Lambda_U^{(I_1,I_2)}(x,y)/
\partial x$ are defined and non decreasing almost everywhere on
$[0,\infty)$.
\end{pro}
\dem   The function $\Lambda_U^{(I_1,I_2)}(x,y)$  is 2-increasing
since a bivariate d.f. is 2-increasing. By Remark \ref{rground} we
conclude that $\Lambda_U^{(I_1,I_2)}(x,y)$ is grounded. Hence,
applying Lemma 2.1.5. in Nelsen (\cite{nelsen_cop}, 2006) we have,
for $(x,y),(x^*,y^*)\in[0,\infty]^2\backslash\{(\infty,\infty)\}$,
$$
\begin{array}{rl}
&|\Lambda_U^{(I_1,I_2)}(x,y)-\Lambda_U^{(I_1,I_2)}(x^*,y^*)|\vspace{0.35cm}\\
\leq&{\dst\lim_{t\to\infty}}t\Big(|P(M(I_1)>1-\frac{x}{t})-P(M(I_1)>1-\frac{x^*}{t})|+|P(M(I_2)>1-\frac{y}{t})
-P(M(I_2)>1-\frac{y^*}{t})|\Big) \vspace{0.35cm}\\
\leq&|I_1||x-x^*|+|I_2||y-y^*|.
\end{array}
$$
Now, the proof is straightforward from Theorem 3 in Schmidt and
Stadtm\"{u}ller (\cite{schmidt+stadt}, 2006). \fdem\\

Remark \ref{rground} and Propositions \ref{pbounds}.(i) and
\ref{pderiv} extend, respectively,  Theorems 1.i), 2.i) and 3 of
Schmidt and Stadtm\"{u}ller (\cite{schmidt+stadt}, 2006). Moreover,
given the above mentioned relation between
$\Lambda_U^{(I_1,I_2)}(x,y)$ and the bivariate upper-tail dependence
function $\Lambda^S_{(F_1(X_1),F_2(X_2))}(x,y)$ in (\ref{utdfS}),
the properties ii)-v) of Theorems 1 and 2 of Schmidt and
Stadtm\"{u}ller (\cite{schmidt+stadt}, 2006) are straightforward for
$\Lambda_U^{(I_1,I_2)}(x,y)$.

\vspace{0.95cm}


We now discuss the case of tail independence between
$M(I_1)$ and $M(I_2)$ and hence extend our context beyond a MEV distribution. \\

Notice that, in case of tail dependence between r.v.'s $F_1(X_1)$
and $F_2(X_2)$, the mapping
\begin{eqnarray}\label{map}
t\mapsto
P\Big(F_1(X_1)>1-\frac{x}{t},F_2(X_2)>1-\frac{y}{t}\Big)
\end{eqnarray}
 is regularly varying of order $-1$ at $\infty$, and so an
homogeneity property holds for large $t$. However, if
$(F_1(X_1),F_2(X_2))$ is tail independent, this latter does not hold
and an adjusted homogeneity property can be obtained by assuming
that
 (\ref{map})
is regularly varying of order $-1/\eta$ at $\infty$, $\eta<1$ (the
case $\eta=1$ corresponds to tail dependence). Coefficient $\eta$ is
the \emph{coefficient of tail dependence} introduced in Ledford and
Tawn (1996, 1997).

Thus being, if we assume that (\ref{map}) is regularly varying of
order $-1/\eta$ at $\infty$, i.e.,
\begin{eqnarray}\label{regvar}
\lim_{t\to \infty}
\frac{P\big(F_1(X_1)>1-x/t,F_2(X_2)>1-y/t\big)}{P\big(F_1(X_1)>1-1/t,F_2(X_2)>1-1/t\big)}=c^*(x,y)
\end{eqnarray}
for $(x,y)\in[0,\infty)^2$, where $c^*$ is homogeneous of order
$1/\eta$ for some $\eta\in (0,1]$ and $c^*(1,1)=1$, then $t\mapsto
P\big(F_1(X_1)>1-1/t,F_2(X_2)>1-1/t\big)$ is regularly varying at
$\infty$ with index $-1/\eta$ (choose $x=y$ in (\ref{regvar})), and
hence we can write
\begin{eqnarray}\label{eta}
P\big(F_1(X_1)>1-1/t,F_2(X_2)>1-1/t\big)=t^{-1/\eta}L(t)
\end{eqnarray}
where $L$ is a slowly varying function at $\infty$ (i.e.,
$L(tx)/L(t)\to 1$, as $t\to \infty$, for any $x>0$). Observe that
$\eta$ dominates the speed of convergence of
$P\big(F_1(X_1)>1-1/t,F_2(X_2)>1-1/t\big)$ to $0$. If $\eta<1$ then
$F_1(X_1)$ and $F_2(X_2)$ (and thus $X_1$ and $X_2$) are
asymptotically independent (or tail independent). In this case, the
tail dependence coefficient $\lambda$ in (\ref{tdc}) is null.
Conversely, asymptotic dependence holds if $\eta=1$ and $L(t)\to
a>0$, as $t\to \infty$, and we have $\lambda>0$. If $\eta=1/2$ we
have (almost) independence (perfect independence if $L(t)=1$ and
(\ref{regvar}) holds with $c^*(x,y)=xy$). The cases $\eta\in(0,1/2)$
and $\eta\in (1/2,1)$ correspond to asymptotically negative
independence and to asymptotically positive independence,
respectively. Roughly speaking, coefficient $\eta$ governs a kind of
a pre-asymptotic tail behavior that allows to better estimate the
probability of extreme events in case of tail independence. A
bivariate extreme value distribution (BEV) allows only tail
dependence ($\eta=1$) or independence ($\eta=1/2$), since
$$
P(F_1(X_1)>1-1/t,F_2(X_2)>1-1/t)\dst\sim (2-l^{(\{1\},\{2\})}(1,1))/t+((l^{(\{1\},\{2\})}(1,1))^2/2-1)/t^2
$$
as $t\to\infty$. For a discussion on this topic see, for instance,
Ledford and Tawn (\cite{led+tawn1}, 1996), Draisma \emph{et al.}
(\cite{draisma+}, 2004) and Drees and M\"{u}ller (\cite{drees},
2008).


\vspace{0.25cm}

\nid Now assume that (\ref{regvar}) holds for random pair
$(M(I_1),M(I_2))$, i.e.,
\begin{eqnarray}\label{regvarM}
\lim_{t\to \infty}
\frac{P\big(M(I_1)>1-x/t,M(I_2)>1-y/t\big)}{P\big(M(I_1)>1-1/t,M(I_2)>1-1/t\big)}=c_{(I_1,I_2)}(x,y)
\end{eqnarray}
for $(x,y)\in[0,\infty)^2$, where $c_{(I_1,I_2)}$ is homogeneous of
order $1/\eta_{(I_1,I_2)}$ for some $\eta_{(I_1,I_2)}\in (0,1]$ and
$c_{(I_1,I_2)}(1,1)=1$. Taking $x=y$ in (\ref{regvarM}), one obtains
that $P\big(M(I_1)>1-1/t,M(I_2)>1-1/t\big)$ is regularly varying at
$\infty$, i.e.,
\begin{eqnarray}\label{etaxy}
P\big(M(I_1)>1-1/t,M(I_2)>1-1/t\big)=t^{-1/\eta_{(I_1,I_2)}}L_{(I_1,I_2)}(t),
\end{eqnarray}
where $L_{(I_1,I_2)}(t)$ is a slowly varying function at $\infty$.
Coefficient $\eta_{(I_1,I_2)}$ is now a measure of the speed of
convergence of $P\big(M(I_1)>1-1/t,M(I_2)>1-1/t\big)$ to $0$ and is,
therefore, a coefficient of tail dependence between $M(I_1)$ and
$M(I_2)$, with analogous conclusions derived for $\eta$ above.
Similarly, in a MEV we obtain, as $t\to\infty$,
$$
P(M(I_1>1-1/t,M(I_2)>1-1/t)\sim (\epsilon_{I_1}+\epsilon_{I_2}-\epsilon_{I_1\cup I_2})/t+(\epsilon_{I_1\cup I_2}^2-
\epsilon_{I_1}^2-\epsilon_{I_2}^2)/(2t^2).
$$
Hence it only occurs asymptotic dependence whenever
$\epsilon_{(I_1,I_2)}=\epsilon_{I_1}+\epsilon_{I_2}-\epsilon_{I_1\cup
I_2}>0$ (with $\eta_{(I_1,I_2)}=1$), and otherwise independence ($\eta_{(I_1,I_2)}=1/2$).\\



\nid In the next result we compute $\eta_{(I_1,I_2)}$ and found that
it is given by the maximum coefficient $\eta_{\{i\},\{j\}}$,
$\forall i\in I_1,j\in I_2$.

\vspace{0.25cm}

\begin{pro}\label{p.3}
Suppose that (\ref{etaxy}) holds and
\begin{eqnarray}\label{etamultivar}
P\Big(\min_{i\in I,j\in
J}(F_i(X_i),F_j(X_j))>1-1/t\Big)=t^{-1/\eta_{I,J}}L_{\eta_{I,J}}(t)
\end{eqnarray}
holds for all $\emptyset\not=I\subset I_1$ and
$\emptyset\not=J\subset I_2$, where $L_{\eta_{I,J}}$ is a slowly
varying function at $\infty$.
Then $\eta_{(I_1,I_2)}=\max\{\eta_{\{i\},\{j\}}:i\in I_1,j\in
I_2\}$.
\end{pro}
\dem First observe that if $I'\subset I$ and $J'\subset J$ then
\begin{eqnarray}\label{desigetamultivar}
1\geq t^{-1/\eta_{I',J'}}L_{\eta_{I',J'}}(t)\geq t^{-1/\eta_{I,J}}L_{\eta_{I,J}}(t).
\end{eqnarray}

We have that
\begin{eqnarray}\label{p.3.1}
\begin{array}{rl}
&P(\bigvee_{i\in I_1}F_i(X_i)>1-1/t,\bigvee_{j\in I_2}F_j(X_j)>1-1/t)\vspace{0.35cm}\\
=&
P(\bigcup_{i\in I_1}\{F_i(X_i)>1-1/t,\bigcup_{j\in I_2}\{F_j(X_j)>1-1/t\}\})\vspace{0.35cm}\\
=&\sum_{\emptyset\not =S\subseteq I_1}(-1)^{|S|+1}P(\bigcap_{i\in S}\{F_i(X_i)>1-1/t\},\bigcup_{j\in I_2}\{F_j(X_j)>1-1/t\}) \vspace{0.35cm}\\
=& \sum_{\emptyset\not =S\subseteq I_1}\sum_{\emptyset\not =T\subseteq I_2}(-1)^{|S|+|T|}
P(\bigcap_{i\in S}\{F_i(X_i)>1-1/t\},\bigcap_{j\in T}\{F_j(X_j)>1-1/t\}) \vspace{0.35cm}\\
=&        \sum_{\emptyset\not =S\subseteq I_1}\sum_{\emptyset\not =T\subseteq I_2}(-1)^{|S|+|T|}
t^{-1/\eta_{S,T}}L_{\eta_{S,T}}(t),
\end{array}
\end{eqnarray}
where in the last equality we have applied (\ref{etamultivar}). Let
\begin{eqnarray}\label{etas0t0}
\eta=\max_{\emptyset\not =S\subseteq I_1 \emptyset\not
=T\subseteq I_2}\eta_{S,T}
\end{eqnarray}

From (\ref{p.3.1}) and (\ref{etas0t0}) we have that
\begin{eqnarray}\label{p.3.2}
\begin{array}{rl}
&P(\bigvee_{i\in I_1}F_i(X_i)>1-1/t,\bigvee_{j\in I_2}F_j(X_j)>1-1/t)\vspace{0.35cm}\\
=&  t^{-1/\eta}L_{\eta}(t)    \sum_{\emptyset\not =S\subseteq I_1}\sum_{\emptyset\not =T\subseteq I_2}(-1)^{|S|+|T|} A_{S,T}(t)
\end{array}
\end{eqnarray}
where $A_{S,T}(t)=t^{-(1/\eta_{S,T}-1/\eta)}L^*_{\eta_{S,T}}(t)$ and
$L^*_{\eta_{S,T}}(t)=L_{\eta_{S,T}}(t)/L_{\eta}(t)$ is a slowly
varying function. Observe that, if $S'\subset S$ and $T'\subset T$,
then $+\infty>A_{S',T'}(t)\geq A_{S,T}(t)$ and, by the definition of
$\eta$, we have $A_{S,T}(t)=1$ or $A_{S,T}(t)\to 0$ as $t\to
\infty$, for all $S\subset I_1$ and $T\subset I_2$. Therefore,
$$
P(M(I_1)>1-1/t,M(I_2)>1-1/t)\sim t^{-1/\eta}L_{\eta}(t).
$$
Moreover, considering $\eta=\eta_{S_0,T_0}$ for some $S_0\subset
I_1,T_0\subset I_2$, and so $A_{S_0,T_0}(t)=1\leq
A_{\{i\},\{j\}}(t)$, $\forall\; i\in S_0,j\in T_0$, we must have
$A_{\{i\},\{j\}}(t)=1$, $\forall\; i\in S_0,j\in T_0$. Then $\eta
=\eta_{\{i\},\{j\}}$, $\forall\; i\in S_0,j\in T_0$ and $\eta\leq
\max_{i\in I_1,j\in I_2}\eta_{\{i\},\{j\}}$. But, by
(\ref{etas0t0}), $\eta\geq \max_{i\in I_1,j\in
I_2}\eta_{\{i\},\{j\}}$ which leads to
the result. \fdem \\

In the following we present some examples where tail independence
takes place. 

\begin{ex}\label{ex3}

Consider $\{V_n\}_{n\geq 1}$ an i.i.d.\hs sequence of r.v.'s with
distribution $\mathtt{U}(0,1)$ and $\mathbf{X}=(X_1,X_2,X_3,X_4)$ a
random vector such that, $X_1=\min(V_3,V_2,V_1)$,
$X_2=\min(V_4,V_2,V_1)$, $X_3=\min(V_4,V_3,V_1)$ and $X_4=V_5$.
Observe that, for $0\leq x\leq 1$,
$F_{X_1}(x)=1-(1-x)^3=F_{X_2}(x)=F_{X_3}(x)$ and $F_{X_4}(x)=x$ and
hence
$F^{-1}_{X_1}(x)=1-(1-x)^{1/3}=F^{-1}_{X_2}(x)=F^{-1}_{X_3}(x)$ and
$F^{-1}_{X_4}(x)=x$. Consider $I_1=\{1,2\}$ and $I_2=\{3,4\}$.

We have successively,
$$
\begin{array}{l}
    P(F_1(X_1)>1-t^{-1},F_3(X_3)>1-t^{-1})=P(F_2(X_2)>1-t^{-1},F_3(X_3)>1-t^{-1})=   t^{-4/3} ,
\end{array}
$$
and
$$
\begin{array}{l}
    P(F_1(X_1)>1-t^{-1},F_4(X_4)>1-t^{-1})=P(F_2(X_2)>1-t^{-1},F_4(X_4)>1-t^{-1})=   t^{-2}    .
\end{array}
$$

Hence, by Proposition \ref{p.3}, we must derive $\eta_{(\{1,2\},\{3,4\})}=3/4$.\\

In fact, applying (\ref{p.3.1}), after some calculations we have
$$
\begin{array}{rl}
&P(M(I_1)>1-t^{-1}x,M(I_2)>1-t^{-1}y)\vspace{0.35cm}\\
=&
\left\{\begin{array}{ll}
2t^{-4/3}xy^{1/3}+2t^{-2}xy-t^{-4/3}x^{4/3}-2t^{-7/3}xy^{4/3}-2t^{-7/3}x^{4/3}y&,\,x\leq y\vspace{0.15cm}\\
t^{-4/3}yx^{1/3}+2t^{-2}xy-3t^{-7/3}x^{1/3}y^{2}-t^{-7/3}x^{4/3}y&,\,x> y.
\end{array}\right.
\end{array}
$$
According to (\ref{etaxy}), coefficient $\eta_{(I_1,I_2)}$ can be
obtained by taking $x=y=1$ in the expression above, and by
(\ref{regvarM}) we obtain
$$
c_{(\{1,2\},\{3,4\})}(x,y)=\left\{\begin{array}{ll}
2xy^{1/3}-x^{4/3}&,\,x\leq y\vspace{0.15cm}\\
yx^{1/3}&,\,x> y.
\end{array}\right.
$$
which is homogeneous of order $4/3$.

Similarly, if we consider $I_1=\{1,2,3\}$ and $I_2=\{4\}$ we obtain
$\eta_{(\{1,2,3\},\{4\})}=1/2$ and $c_{(\{1,2,3\},\{4\})}(x,y)=xy$,
and if $I_1=\{1\}$ and $I_2=\{2,3,4\}$ we have
$\eta_{(\{1\},\{2,3,4\})}=3/4$ and
$$
c_{(\{1\},\{2,3,4\})}(x,y)=\left\{\begin{array}{ll}
xy^{1/3}&,\,x\leq y\vspace{0.15cm}\\
2yx^{1/3}-y^{4/3}&,\,x> y
\end{array}\right.=c_{(\{1,2\},\{3,4\})}(y,x).\\
$$\\
%
%
\end{ex}

\begin{ex}\label{ex4}

Consider $\mathbf{X}=(X_1,...X_d)$ a standard $d$-variate Gaussian
random vector with positive definite correlation matrix.
The bivariate tail-dependence structure is given by
\begin{eqnarray}\label{}
P(F_i(X_i)>1-1/t,F_j(X_j)>1-1/t)\sim C_{\rho_{i,j}} t^{-2/(1+\rho_{i,j})}(\log(t))^{-\rho_{i,j}/(1+\rho_{i,j})},\textrm{ as $t\to \infty$,}
\end{eqnarray}
for $i,j\in\{1,...,d\}$, $i<j$, where
$\rho_{i,j}=corr(X_i,X_j)\not\in\{-1,1\}$ and
$$C_{\rho_{i,j}}=(1+\rho_{i,j})^{3/2}(1-\rho_{i,j})^{-1/2}(4\pi)^{-\rho_{i,j}/(1+\rho_{i,j})}.$$
Hence (\ref{etamultivar}) holds for $I=\{i\}$ and $J=\{j\}$ with
$\eta_{i,j}=(1+\rho_{i,j})/2$ (see Ledford and Tawn
\cite{led+tawn1}, 1996; Draisma \emph{et al.} \cite{draisma+},
2004). According to Hua and Joe (\cite{hua+joe}, 2011),
(\ref{etamultivar}) also holds for non-empty sets
$I_1,I_2\subset\{1,...,d\}$. If we consider
$\rho_{(I_1,I_2)}=\max\{\rho_{i,j}:i\in I_1,j\in I_2\}$ then, by
Proposition \ref{p.3}, we find
$\eta_{(I_1,I_2)}=(1+\rho_{(I_1,I_2)})/2$, provided the left-hand
side of (\ref{etaxy}) is non-null.
\end{ex}

\section{Estimation}\label{sestim}

Several estimators for the bivariate stable tail dependence function
in (\ref{stf}) or even for the more general $d$-variate stable tail
dependence function
\begin{eqnarray}\label{stfmulti}
\lim_{t\to\infty}tP\Big(F_1(X_1)>1-\frac{x_1}{t} \vee ... \vee F_d(X_d)>1-\frac{x_d}{t}\Big)
\end{eqnarray}
have been considered in literature. For a survey, see Krajina (2010)
\cite{krajina}.  According to relation (\ref{stdfnossa}), they can
be applied to our function $l^{(I_1,I_2)}(x^{-1},y^{-1})$.

We remark that these are based on asymptotic results that depend on
a sequence of positive integers, $\{k_n\}$, going to infinity at a
lower rate than $n$. For instance, the estimator based on
(\ref{stf}) by plugging-in the respective empirical counterparts
given by
$$
\frac{n}{k_n}P_n\Big(\widehat{F}_1(X_1)>1-\frac{k_n}{n}x \vee \widehat{F}_2(X_2)>1-\frac{k_n}{n}y\Big)
=\frac{1}{k_n}\sum_{i=1}^n\mathbf{1}_{\{\widehat{F}_1(X_1)>1-\frac{k_n}{n}x \vee \widehat{F}_2(X_2)>1-\frac{k_n}{n}y\}},
$$
where $\widehat{F}_l(u)=n^{-1}\sum_{k=1}^n\mathbf{1}_{\{X_k\leq
u\}}$ is the empirical d.f. of $F_l$, $l=1,2$, is consistent and
asymptotically normal if  $\{k_n\}$ is an intermediate sequence,
i.e., $k_n\to\infty$ and $k_n/n\to 0$, as $n\to\infty$ (Huang 1992
\cite{huang}).
 The choose of the value $k$ in the sequence $\{k_n\}$ that allows
the better trade-off between bias and variance is of major
difficulty, since small values of $k$ come along with a large
variance whenever an increasing $k$ results in a strong bias.
Therefore, simulation studies have been carried out in order to find
the best value of $k$
that allows this compromise.\\


As mentioned before, the upper-tail dependence function in
(\ref{utdf}) can be viewed as an extension of the bivariate
upper-tail dependence function of Schmidt and Stadtm\"{u}ller (2006)
given in (\ref{utdfS}),
by taking in this limit the random pair $(M(I_1),M(I_2))$ instead of $(F(X_1),F(X_2))$. 
The estimators considered in  Schmidt and Stadtm\"{u}ller (2006),
for which strong consistency and asymptotic normality have been
established, allow to estimate our function
$\Lambda_{U}^{(I_1,I_2)}(x,y)$, as well as coefficient
$\epsilon_{(I_1,I_2)}=\Lambda_{U}^{(I_1,I_2)}(1,1)$. However they
are also based on asymptotic results with the same
drawback of including an intermediate sequence, already referred above.\\

In order to overcome this problem, we shall present a totally
different and very simple approach. More precisely, the following
result suggests an estimation procedure for the $d$-variate stable
tail dependence function in (\ref{stfmulti}) that only evolves a
sample mean.


\begin{pro}\label{p.2}
Under the conditions of Proposition \ref{p.1}, we have, for
$l(x_1,...,x_d)=-\log F(x_1,...,x_d)$,
\begin{eqnarray}\nn
l(x_1,...,x_d)=\frac{E(F_1(X_1)^{x_1}\vee ...\vee F_d(X_d)^{x_d})}{1-E(F_1(X_1)^{x_1}\vee ...\vee F_d(X_d)^{x_d})}.
\end{eqnarray}
\end{pro}
\dem 
Consider  for $G(x)=\exp(-1/x)$. Observe that
 $$
E(G(x_1X_1)\vee ... \vee G(x_dX_d))=E(G( x_1X_1\vee ... \vee x_dX_d))
 $$
and the d.f. of $x_1X_1\vee ... \vee x_dX_d$  is given by
\begin{eqnarray}\label{pstdfmulti2}
\begin{array}{rl}
 P( x_1X_1\vee ... \vee x_dX_d\leq u)=&P(X_1\leq u/x_1, ...,X_d\leq u/x_d)\vspace{0.35cm}\\
 = &F\big(ux_1^{-1},...,ux_d^{-1}\big)\vspace{0.35cm}\\
 = & \exp \big(-l\big(ux_1^{-1},...,ux_d^{-1}\big) \big).
\end{array}
\end{eqnarray}
Hence
\begin{eqnarray}\label{pstdfmulti3}
\begin{array}{rl}
 &E(G( x_1X_1\vee ... \vee x_dX_d))\vspace{0.35cm}\\
 =&\dst\int_0^\infty \exp(-u^{-1})\exp \big(-l\big(ux_1^{-1},...,ux_d^{-1}\big) \big)\frac{d}{du}
  \big(-l\big(ux_1^{-1},...,ux_d^{-1}\big) \big)
 \vspace{0.35cm}\\
 =&\dst\int_0^\infty \exp(-u^{-1})\exp \big(-u^{-1}l\big(x_1^{-1},...,x_d^{-1}\big) \big)\frac{d}{du}
  \big(-u^{-1}l\big(x_1^{-1},...,x_d^{-1}\big) \big)\vspace{0.35cm}\\
 = & l\big(x_1^{-1},...,x_d^{-1}\big)\dst\int_0^\infty \exp \big(-u^{-1}\big(1+l\big(x_1^{-1},...,x_d^{-1}\big)\big) \big)
  u^{-2} du\vspace{0.35cm}\\
=& \dst\frac{l\big(x_1^{-1},...,x_d^{-1}\big)}{1+l\big(x_1^{-1},...,x_d^{-1}\big)}.
\end{array}
\end{eqnarray}
Now just observe that $G(  x_1X_1\vee ... \vee x_dX_d)\dst\mathop{=}^dF_1(X_1)^{1/x_1}\vee...\vee F_d(X_d)^{1/x_d}$. \fdem\\

\begin{rem}
Observe that the $d$-variate stable tail dependence function in
(\ref{stfmulti}) corresponds to $-\log F(x_1^{-1},...,x_d^{-1})$.
\\
\end{rem}
By applying Proposition \ref{p.2} with $x_j$ replaced by $x_j^{-1}$,
$j=1,...,d$, we get the following corollary.

\begin{cor}\label{c.1}
Under the conditions of Proposition \ref{p.1}, we have
\begin{eqnarray}\label{xepsilon1}
x\epsilon_{I_1}\equiv l^{(I_1,\emptyset)}(x^{-1},x^{-1})=\frac{E(M(I_1)^{1/x})}{1-E(M(I_1)^{1/x})},
\end{eqnarray}
\begin{eqnarray}\label{yepsilon2}
y\epsilon_{I_2}\equiv l^{(\emptyset,I_2)}(y^{-1},y^{-1})=\frac{E(M(I_2)^{1/y})}{1-E(M(I_2)^{1/y})}
\end{eqnarray}
and
\begin{eqnarray}\label{pvalmed}
l^{(I_1,I_2)}(x^{-1},y^{-1})=\frac{E(M(I_1)^{1/x}\vee M(I_2)^{1/y})}{1-E(M(I_1)^{1/x}\vee M(I_2)^{1/y})}.
\end{eqnarray}\\
\end{cor}

Consider  the estimators derived from Proposition \ref{p.1}
 and Corollary \ref{c.1} by plugging-in the respective sample means,
 respectively,
\begin{eqnarray}\label{estimadores2}
\begin{array}{c}
\widetilde{l}(x_1,...,x_d)=
\dst\frac{\overline{F_1(X_1)^{x_1}\vee ...\vee F_d(X_d)^{x_d}}}{1-\overline{F_1(X_1)^{x_1}\vee ...\vee F_d(X_d)^{x_d}}},
\end{array}
\end{eqnarray}
and
\begin{eqnarray}\label{estimadores1}
\begin{array}{ccc}
\dst x\widetilde{\epsilon}_{I_1}=\frac{\overline{M(I_1)^{1/x}}}{1-\overline{M(I_1)^{1/x}}},&\hspace{-0.125cm}\dst
 y\widetilde{\epsilon}_{I_2}=\frac{\overline{M(I_2)^{1/y}}}{1-\overline{M(I_2)^{1/y}}}\,\textrm{ and}&\hspace{-0.125cm} \dst
\widetilde{l}^{(I_1,I_2)}(x^{-1},y^{-1})=\frac{\overline{M(I_1)^{1/x}\vee M(I_2)^{1/y}}}{1-\overline{M(I_1)^{1/x}\vee M(I_2)^{1/y}}}.
\end{array}
\end{eqnarray}
where
\begin{eqnarray}\label{estimconhec}
\begin{array}{cc}
\dst\overline{M(I_1)^{1/x}}=\frac{1}{n}\sum_{i=1}^n\bigvee_{j\in I_1}F_j(X^{(i)}_j)^{1/x},&
\overline{M(I_2)^{1/y}}=\dst\frac{1}{n}\sum_{i=1}^n\bigvee_{j\in I_2}F_j(X^{(i)}_j)^{1/y}
\end{array}
\end{eqnarray}
and
\begin{eqnarray}\label{estimconhec2}
\begin{array}{c}
\overline{M(I_1)^{1/x}\vee M(I_2)^{1/y}}=\dst\frac{1}{n}\sum_{i=1}^n\Big(\bigvee_{j\in I_1}F_j(X^{(i)}_j)^{1/x}\vee
\bigvee_{j\in I_2}F_j(X^{(i)}_j)^{1/y} \Big).
\end{array}
\end{eqnarray}
We will consider two situations: the first one for known margins and the second one for unknown margins.  \\

In case the margins are known, they become unit Fr\'{e}chet by
transformation $-1/\log F_j(X_j)$ for $j\in I\subset\{1,...,d\}$.

It is quite straightforward to deduce the consistency and asymptotic
normality of estimators (\ref{estimadores2}) and
(\ref{estimadores1}) by the well-known Delta Method.
\begin{pro}\label{pestim}
Under the conditions of Proposition \ref{p.1}, we have
\begin{eqnarray}\label{normxepsilon1}
\sqrt{n}(\widetilde{l}(x_1,...,x_d)-l(x_1,...,x_d))\to
N(0,\sigma^2),
\end{eqnarray}
where $\widetilde{l}(x_1,...,x_d)$ is the estimator derived from
Proposition \ref{p.2} by plugging-in the respective sample mean
given in (\ref{estimadores2}) and
\begin{eqnarray}\nn
\begin{array}{c}
\sigma^2
=\frac{l(x_1,...,x_d)\big(1+l(x_1,...,x_d)\big)^2}{\big(2+l(x_1,...,x_d)\big)}.
\end{array}
\end{eqnarray}
\end{pro}
\dem  Let $Y_i$, $i=1,...,n$, be independent copies of
$Y=F_1(X_1)^{x_1}\vee ...\vee F_d(X_d)^{x_d}$. We have that
$\sqrt{n}(\overline{Y}-\mu_Y)\to N(0,\sigma^2_Y)$, where
$\mu_Y=E(F_1(X_1)^{x_1}\vee ...\vee F_d(X_d)^{x_d})$ and
$\sigma^2_Y=Var(F_1(X_1)^{x_1}\vee ...\vee F_d(X_d)^{x_d})$. By a
similar reasoning of (\ref{pstdfmulti3}) we derive
$$
\begin{array}{c}
E((F_1(X_1)^{x_1}\vee ...\vee F_d(X_d)^{x_d})^2)=\frac{l(x_1,...,x_d)}{2+l(x_1,...,x_d)}
\end{array}
$$
and hence, 
$$
\begin{array}{c}
Var((F_1(X_1)^{x_1}\vee ...\vee F_d(X_d)^{x_d})^2)=\frac{l(x_1,...,x_d)}
{\big(2+l(x_1,...,x_d)\big)\big(1+l(x_1,...,x_d)\big)^2}.
\end{array}
$$

\nid Let $g(x)=(1-x)^{-1}-1$. We have $[g'(\mu_Y)]^2=(1-\mu_Y)^{-4}$
and, by the Delta Method,
$\sqrt{n}(g(\overline{Y})-x\epsilon_{I_1})\to
N(0,\sigma^2_Y(1-\mu_Y)^{-4})$. \fdem\\

\begin{cor}\label{cestim}
Under the conditions of Proposition \ref{p.1}, we have
\begin{eqnarray}\label{normxepsilon1}
\sqrt{n}(x\widetilde{\epsilon}_{I_1}-x\epsilon_{I_1})\to
N(0,\sigma^2_1),
\end{eqnarray}
\begin{eqnarray}\label{normyepsilon2}
\sqrt{n}(y\widetilde{\epsilon}_{I_2}-y\epsilon_{I_2})\to
N(0,\sigma^2_2)
\end{eqnarray}
and
\begin{eqnarray}\label{normpvalmed}
\sqrt{n}(\widetilde{l}^{(I_1,I_2)}(x^{-1},y^{-1})-l^{(I_1,I_2)}(x^{-1},y^{-1}))\to
N(0,\sigma^2_3),
\end{eqnarray}
where $x\widetilde{\epsilon}_{I_1}$, $y\widetilde{\epsilon}_{I_2}$
and $\widetilde{l}^{(I_1,I_2)}(x^{-1},y^{-1})$ are given in
(\ref{estimadores1}) and
\begin{eqnarray}\label{varxepsilon1}
\sigma^2_{1}=\frac{x\epsilon_{I_1}\big(1+x\epsilon_{I_1}\big)^2}
{\big(2+x\epsilon_{I_1}\big)},
\end{eqnarray}
\begin{eqnarray}\label{varyepsilon2}
\sigma^2_{2}=\frac{y\epsilon_{I_2}\big(1+y\epsilon_{I_2}\big)^2}
{\big(2+y\epsilon_{I_2}\big)}
\end{eqnarray}
and
\begin{eqnarray}\label{varpvalmed}
\sigma^2_{3}=\frac{l^{(I_1,I_2)}\big(x^{-1},y^{-1}\big)\big(1+l^{(I_1,I_2)}\big(x^{-1},y^{-1}\big)\big)^2}
{\big(2+l^{(I_1,I_2)}\big(x^{-1},y^{-1}\big)\big)}.
\end{eqnarray}

\end{cor}

%

Based on the definition in (\ref{utdf}), a natural estimator for the
upper-tail dependence function is
\begin{eqnarray}\label{estimutdf}
\widetilde{\Lambda}_{U}^{(I_1,I_2)}(x,y)=x\widetilde{\epsilon_{I_1}}+y\widetilde{\epsilon_{I_2}}-\widetilde{l}^{(I_1,I_2)}(x^{-1},y^{-1}),
\end{eqnarray}
with $x\widetilde{\epsilon_{I_1}}$, $y\widetilde{\epsilon_{I_2}}$
and $\widetilde{l}^{(I_1,I_2)}(x^{-1},y^{-1})$ stated in
(\ref{estimadores1}). Hence we have the following estimator for the
extremal coefficient of dependence between
$\mathbf{X}_{I_1}$ and $\mathbf{X}_{I_2}$:
\begin{eqnarray}\label{estimextcoef}
\widetilde{\epsilon}_{(I_1,I_2)}=\widetilde{\epsilon_{I_1}}+\widetilde{\epsilon_{I_2}}-\widetilde{\epsilon}_{I_1\cup I_2}.
\end{eqnarray}
where $\widetilde{\epsilon}_{I_1\cup
I_2}=\widetilde{l}^{(I_1,I_2)}(1,1)$.

\begin{pro}\label{strongconsistentcy}
Estimators $\widetilde{l}(x_1,...,x_d)$ and
$\widetilde{\Lambda}_{U}^{(I_1,I_2)}(x,y)$ in (\ref{estimadores2})
 and (\ref{estimutdf}), respectively, are strong consistent.
Consequently, the same holds for $\widetilde{\epsilon}_{(I_1,I_2)}$
in (\ref{estimextcoef}).
\end{pro}
\dem Just observe that, as the sample mean $\overline{M(I_1)^{1/x}}$
converges almost surely to the mean value $E(M(I_1)^{1/x})$, i.e.,
$\overline{M(I_1)^{1/x}}\dst\mathop{\longrightarrow}^{a.s.}
E(M(I_1)^{1/x})$, then
$x\widetilde{\epsilon_{I_1}}=g(\overline{M(I_1)^{1/x}})\dst\mathop{\longrightarrow}^{a.s.}
x\epsilon_{I_1}=g(E(M(I_1)^{1/x}))$, where $g(x)=(1-x)^{-1}-1$.
Analogously for $y\widetilde{\epsilon_{I_2}}$,
$\widetilde{l}^{(I_1,I_2)}(x^{-1},y^{-1})$ and
$\widetilde{l}(x_1,...,x_d)$. Now, the  strong consistency of
$\widetilde{\Lambda}_{U}^{(I_1,I_2)}(x,y)$ is straightforward from
$$
\begin{array}{rl}
  |\widetilde{\Lambda}_{U}^{(I_1,I_2)}(x,y)-\Lambda_{U}^{(I_1,I_2)}(x,y) |
  \leq |x\widetilde{\epsilon_{I_1}}-x\epsilon_{I_1}|+ |y\widetilde{\epsilon_{I_2}}-y\epsilon_{I_2}|
  + |\widetilde{l}^{(I_1,I_2)}(x^{-1},y^{-1})-l^{(I_1,I_2)}(x^{-1},y^{-1})|.\,\,$\fdem$
\end{array}
$$
\\

  Now consider $\widehat{F}_j$ the empirical d.f. of $F_j$,
$j=1,...,d$,
$$
\widehat{F}_j(u)=\frac{1}{n+1}\sum_{k=1}^n\mathbf{1}_{\{X_j^{(k)}\leq
u\}}.
$$
The denominator $n+1$ instead of $n$ in the empirical d.f. concerns
estimation accuracy and other modifications can be used. For a
discussion see, for instance, Beirlant et al. \cite{beirl+} (2004).

In case of unknown margins, we can replace $F_j$ by the respective
empirical d.f. $\widehat{F}_j$, $j=1,...,d$, in (\ref{estimadores2})
and (\ref{estimadores1}). More precisely, we have
\begin{eqnarray}\label{estimadores2emp}
\begin{array}{c}
\dst\widehat{l}(x_1,...,x_d)=
\frac{\overline{\widehat{F}_1(X_1)^{x_1}\vee ...\vee \widehat{F}_d(X_d)^{x_d}}}{1-\overline{\widehat{F}_1(X_1)^{x_1}\vee ...\vee \widehat{F}_d(X_d)^{x_d}}},
\end{array}
\end{eqnarray}
as well as,
\begin{eqnarray}\label{estimadores1emp}
\begin{array}{ccc}
\dst x\widehat{\epsilon}_{I_1}=\frac{\overline{\widehat{M}(I_1)^{1/x}}}{1-\overline{\widehat{M}(I_1)^{1/x}}},&\hspace{-0.125cm}\dst
 y\widehat{\epsilon}_{I_2}=\frac{\overline{\widehat{M}(I_2)^{1/y}}}{1-\overline{\widehat{M}(I_2)^{1/y}}}\,\textrm{ and}&\hspace{-0.125cm} \dst
\widehat{l}^{(I_1,I_2)}(x^{-1},y^{-1})=\frac{\overline{\widehat{M}(I_1)^{1/x}\vee \widehat{M}(I_2)^{1/y}}}{1-\overline{\widehat{M}(I_1)^{1/x}\vee \widehat{M}(I_2)^{1/y}}}
\end{array}
\end{eqnarray}
where
\begin{eqnarray}\label{estimdescl}
\begin{array}{c}
\dst\overline{\widehat{F}_1(X_1)^{x_1}\vee ...\vee \widehat{F}_d(X_d)^{x_d}}=\frac{1}{n}\sum_{i=1}^n\bigvee_{j\in
\{1,...,d\}}\widehat{F}_j(X^{(i)}_j)^{x_j},
\end{array}
\end{eqnarray}
\begin{eqnarray}\label{estimdesc}
\begin{array}{cc}
\dst\overline{\widehat{M}(I_1)^{1/x}}=\frac{1}{n}\sum_{i=1}^n\bigvee_{j\in I_1}\widehat{F}_j(X^{(i)}_j)^{1/x},&
\overline{\widehat{M}(I_2)^{1/y}}=\dst\frac{1}{n}\sum_{i=1}^n\bigvee_{j\in I_2}\widehat{F}_j(X^{(i)}_j)^{1/y}
\end{array}
\end{eqnarray}
and
\begin{eqnarray}\label{estimdesc2}
\begin{array}{c}
\overline{\widehat{M}(I_1)^{1/x}\vee \widehat{M}(I_2)^{1/y}}=\dst\frac{1}{n}\sum_{i=1}^n\Big(\bigvee_{j\in I_1}\widehat{F}_j(X^{(i)}_j)^{1/x}\vee
\bigvee_{j\in I_2}\widehat{F}_j(X^{(i)}_j)^{1/y} \Big).
\end{array}
\end{eqnarray}

We still have asymptotic normality of estimators in
(\ref{estimdescl})-(\ref{estimdesc2}) from the following result
stated in Fermanian \emph{et al.} (2002, \cite{ferm}, Theorem 6).

\begin{teo}\emph{(Fermanian \emph{et al.} (2002) \cite{ferm}, Theorem
6)} Let $F$ have continuous marginals and let copula $C_F$ in
(\ref{copula}) have continuous partial derivatives. Then
$$
\frac{1}{\sqrt{n}}\sum_{i=1}^n \{J(\widehat{F}_1(X^{(i)}_1),...,\widehat{F}_d(X^{(i)}_d))
-E(J(F_1(X^{(i)}_1),...,F_d(X^{(i)}_d)))\}\to\int_{[0,1]^d}\mathbb{G}(u_1,...,u_d)dJ(u_1,...,u_d)
$$
in distribution in $\ell^\infty([0,1]^d)$, where the limiting
process and $\mathbb{G}$ are centered Gaussian, and
$J:[0,1]^d\to\mathbb{R}$ is of bounded variation, continuous from
above and with discontinuities of the first kind (Neuhaus, 1971
\cite{neuhaus}).
\end{teo}

The asymptotic normality of estimators (\ref{estimadores2emp}) and
(\ref{estimadores1emp}) is now derived from a general version of the
Delta Method as considered in Schmidt and Stadtm\"{u}ller
\cite{schmidt+stadt} (2006; Theorem 13).\\

We also state strong consistency of estimators
$\widehat{l}(x_1,...,x_d)$ in (\ref{estimadores2emp}) and
\begin{eqnarray}\label{estimutdfemp}
\widehat{\Lambda}_{U}^{(I_1,I_2)}(x,y)=x\widehat{\epsilon_{I_1}}+y\widehat{\epsilon_{I_2}}-\widehat{l}^{(I_1,I_2)}(x^{-1},y^{-1}),
\end{eqnarray}
with $x\widehat{\epsilon_{I_1}}$, $y\widehat{\epsilon_{I_2}}$ and
$\widehat{l}^{(I_1,I_2)}(x^{-1},y^{-1})$ given in
(\ref{estimadores1emp}), and hence of estimator
\begin{eqnarray}\label{estimextcoefemp}
\widehat{\epsilon}_{(I_1,I_2)}=\widehat{\epsilon_{I_1}}+\widehat{\epsilon_{I_2}}-\widehat{\epsilon}_{I_1\cup I_2},
\end{eqnarray}
where $\widehat{\epsilon}_{I_1\cup
I_2}=\widehat{l}^{(I_1,I_2)}(1,1)$.

\begin{pro}
Estimators $\widehat{l}(x_1,...,x_d)$ in (\ref{estimadores2emp}) and
$\widehat{\Lambda}_{U}^{(I_1,I_2)}(x,y)$ in (\ref{estimutdfemp}) are
strong consistent. Therefore, the same holds for estimator
$\widehat{\epsilon}_{(I_1,I_2)}$ in (\ref{estimextcoefemp}).
\end{pro}
\dem The proof runs along the same lines as the one of Proposition
\ref{strongconsistentcy}. We only prove the more general case
$\widehat{l}(x_1,...,x_d)\dst\mathop{\longrightarrow}^{a.s.}l(x_1,...,x_d)$.
Observe that
$$
\begin{array}{rl}
&\dst\Big|\frac{1}{n}\sum_{i=1}^n\bigvee_{j\in{\{1,...,d\}}}\widehat{F_j}(X_j^{(i)})^{x_j}
-E\big(\bigvee_{j\in{\{1,...,d\}}}F_j(X_j)^{x_j}\big)\Big|\vspace{0.35cm}\\
\leq &\dst\Big|\frac{1}{n}\sum_{i=1}^n\bigvee_{j\in{\{1,...,d\}}}\widehat{F_j}(X_j^{(i)})^{x_j}-
\frac{1}{n}\sum_{i=1}^n\bigvee_{j\in{\{1,...,d\}}}F_j(X_j^{(i)})^{x_j}\Big|\vspace{0.35cm}\\
&\dst+\Big|\frac{1}{n}\sum_{i=1}^n\bigvee_{j\in{\{1,...,d\}}}F_j(X_j^{(i)})^{x_j}-E\big(\bigvee_{j\in{\{1,...,d\}}}F_j(X_j)^{x_j}\big)\Big|,
\end{array}
$$
where the second term converges \emph{almost surely} to zero by the
\emph{Strong Law of Large Numbers}.\\

For the first term we have, successively,
$$
\begin{array}{l}
|\frac{1}{n}\sum_{i=1}^n\bigvee_{j\in\{1,...,d\}}\widehat{F}_j(X_j^{(i)})^{x_j}-\frac{1}{n}\sum_{i=1}^n\bigvee_{j\in\{1,...,d\}}F_j(X_j^{(i)})^{x_j}|\vspace{0.35cm}\\
\leq\frac{1}{n}\sum_{i=1}^n\bigvee_{j\in\{1,...,d\}}|\widehat{F}_j(X_j^{(i)})^{x_j}-F_j(X_j^{(i)})^{x_j}|\vspace{0.35cm}\\
\leq\frac{1}{n}\sum_{i=1}^n\sum_{j\in\{1,...,d\}}|\widehat{F}_j(X_j^{(i)})^{x_j}-F_j(X_j^{(i)})^{x_j}|,
\end{array}
$$
 which
converges \emph{almost surely} to zero according to Gilat and
Hill (\cite{gilat+hill}, 1992; proof of Theorem 1.1). \fdem 

\section{Application to financial data}\label{saplic}

\begin{figure}[!htb]
\begin{center}
\begin{tabular}{l}
 \includegraphics[width=12.55cm,height=4.4cm]{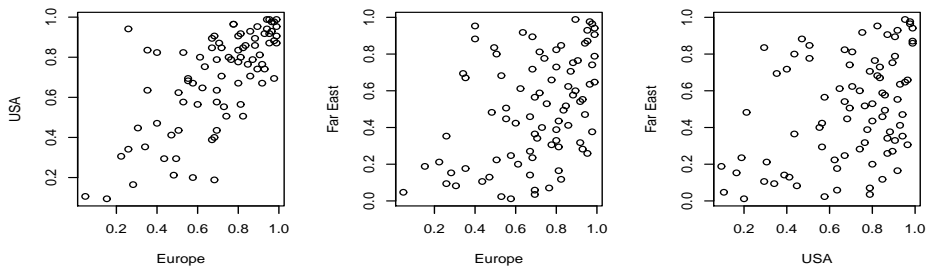}
\end{tabular}
\caption{Scatter plots of the monthly maximums ($84$ data points) in Europe versus USA, Europe versus Far East and USA versus Far East.\label{fig1}}
\end{center}
\end{figure}

In this section we show that tail dependence is present in financial
data. Our analysis is based on negative log-returns of daily closing
values of the stock market indexes, CAC 40 (France), FTSE100 (UK),
SMI (Swiss), XDAX (German), Dow Jones (USA), Nasdaq (USA), SP500
(USA), HSI (China), Nikkei (Japan). The period covered is January
1993 to March 2004. More precisely, we consider the monthly maximums
in each market and group the indexes in Europe (CAC 40, FTSE100,
SMI, XDAX), USA (Dow Jones, Nasdaq) and Far East (HSI, Nikkei). The
scatter plots in Figure \ref{fig1} show the presence of dependence
between the monthly maximums in Europe and USA, Europe and Far East,
USA and Far East, respectively.
We are interested in assessing the amount of tail dependence between
the three big world markets referred: Europe, USA and Far East, and
this can be achieved through the \emph{extremal coefficient of
dependence} $\epsilon_{(I_1,I_2)}$, defined in (\ref{extcoef}). As
we do not know the margins distribution, we use estimator
$\widehat{\epsilon}_{(I_1,I_2)}$ in (\ref{estimextcoefemp}) based on
ranks. In Table \ref{tab1} are the obtained estimates for several
groups, $I_1$ and $I_2$. One can see that the Far East market has
less influence (lower values of the coefficient) but Europe and USA
have a stronger effect on each other and on the respective group of
foreign markets. Observe that the difference between these two
magnitudes of dependence is almost in the proportion 1:2.

\begin{table}[!hbp]
\begin{center}
\begin{tabular}{|c|c|c|}
\hline
$I_1$&$I_2$&$\widehat{\epsilon}_{(I_1,I_2)}$\\
\hline
Europe&USA & 1.008324625
\\
Europe & Far East& 0.568780467
\\
USA & Far East& 0.364358832
\\
Europe & USA $\cup$ Far East& 1.125919957
\\
USA & Europe $\cup$ Far East& 0.921498322
\\
Far East & USA $\cup$ Europe  & 0.481954164
\\
\hline
\end{tabular}
\caption{Estimates of the \emph{extremal coefficient of dependence}
$\widehat{\epsilon}_{(I_1,I_2)}$ for the indicated groups, $I_1$ and
$I_2$. \label{tab1}}
\end{center}
\end{table}

\section{Conclusion}
In this work we introduce a new upper-tail dependence concept for a
random vector which extends the one in Schmidt and Stadtm\"{u}ller
(\cite{schmidt+stadt}, 2006). Our approach weakens the usual imposed
multivariate tail dependence and can be treated with bivariate
techniques. The new function extends the well-known relation of
Huang (\cite{huang} 1992) for a MEV with unit Fréchet marginals, and
gives rise to the so-called \emph{extremal coefficient of
dependence} as it is expressed through the extremal coefficient in
Tiago de Oliveira (\cite{tiago}, 1962-63) and Smith (\cite{smith},
1990). We also enlarge our discussion to tail independence in the
sense of Ledford and Tawn (\cite{led+tawn1,led+tawn2}, 1996, 1997).
At this point we are beyond MEV distributions which only admit tail
dependence or (exact) independence.

In calculating the moments of the r.v.'s involved in our function,
we arrive at very simple estimators whose asymptotic normality is
stated. These can also be applied to the well-known \emph{stable
tail dependence function}. We also prove strong consistency of the
proposed estimators for our measures. We end with an application to
financial data presenting tail dependence.

 \vspace{0.75cm}


\begin{thebibliography}{99}

\bibitem{beirl+} Beirlant, J.,  Goegebeur, Y., Segers, J. e Teugels, J. (2004). \emph{Statistics of Extremes: Theory and Application.} John Wiley.

\bibitem{coles+} Coles, S., Heffernan, J. and Tawn, J. (1999). Dependence measures
for extreme value analysis, Extremes 2: 339-366.

\bibitem{draisma+} Draisma, G.,  Drees, H., Ferreira, A.\hspace{3pt}and de Haan, L. (2004). Bivariate tail estimation: dependence in asymptotic independence. \emph{Bernoulli}, 10, 251-280.


\bibitem{drees} Drees, H., and Müller, P. (2008). Fitting and validation of a
bivariate model for large claims. Insurance: Mathematics and
Economics 42, 638-650.



\bibitem{emb+03} Embrechts, P., Lindskog, F. and McNeil, A. (2003). Modelling
Dependence with Copulas and Applications to Risk Management, In:
Handbook of Heavy Tailed Distibutions in Finance, ed. S. Rachev,
Elsevier, Chapter 8: 329-384.

\bibitem{ferm}Fermanian, J.-D., Radulovi\'{c}, D.,  Wegkamp, M. (2004). Weak
convergence of empirical copula processes. Bernoulli 10(5), 847-860.

\bibitem{gilat+hill} Gilat, D. and Hill, T. (1992) One-sided refinements of the
strong law of large numbers and the Glivenko-Cantelli Theorem, Ann.
Probab. 20 , 1213-1221.

\bibitem{hua+joe} Hua, L., Joe, H. (2004). Tail order and intermediate tail dependence of multivariate copulas.

\bibitem{huang} Huang, X. (1992). Statistics of Bivariate Extreme Values. Ph. D.
thesis, Tinbergen Institute Research Series 22, Erasmus University
Rotterdam.

\bibitem{joe} Joe, H. (1997). Multivariate Models and Dependence Concepts. Chapman
\& Hall, London.

\bibitem{krajina} Krajina, A. (2010). An M-Estimator of Multivariate Tail Dependence.
Tilburg: Tilburg University Press.

\bibitem{led+tawn1}{Ledford, A.\hspace{3pt}and Tawn, J.\hspace{3pt}A.} (1996). Statistics for near independence in multivariate extreme values. {\em Biometrika}, 83, 169-187.

\bibitem{led+tawn2} {Ledford, A.\hspace{3pt}and Tawn, J.\hspace{3pt}A.}
(1997). Modelling Dependence within joint tail regions, J. R. Stat.
Soc. Ser. B Stat. Methodol. 59, 475-499.

\bibitem{li1} Li, H. (2006). Tail dependence of multivariate Pareto distributions,
WSU Mathematics Technical Report 2006-6.
http://www.math.wsu.edu/TRS/.

\bibitem{li2} Li, H. (2008). Tail Dependence
Comparison of Survival Marshall-Olkin Copulas, Methodol. Comput.
Appl. Probab., 10(1), 39-54.

\bibitem{li3} Li, H. (2009). Orthant tail dependence of multivariate extreme value
distributions, J. Multivariate Anal., 100(1), 243-256.

\bibitem{mar+olkin} Marshall, A.W., Olkin, I. (1967). A multivariate exponential
distribution, J. Amer. Statist. Assoc. 62 30-44.

\bibitem{nelsen} Nelsen, R. B. (1996). Nonparametric measures of multivariate
association, in Distribution with fixed marginals and related
topics, IMS Lecture Notes - Monograph Series, vol. 28, 223-232.

\bibitem{nelsen_cop} Nelsen, R.B. (2006). An Introduction to Copulas. Second Edition.
Springer, New York.

\bibitem{neuhaus} Neuhaus, G. (1971) On the weak convergence of stochastic processes
with multidimensional time parameter. Ann. Math. Statist., 42,
1285-1295.

\bibitem{schmid+schmidt} Schmid, F., Schmidt, R. (2007). Multivariate conditional versions of
Spearman's rho and related measures of tail dependence. J.
Multivariate Anal., 98, 1123-1140.

\bibitem{schmidt+stadt} Schmidt, R., Stadtm\"{u}ller, U. (2006). Nonparametric estimation of
tail dependence, The Scandinavian Journal of Statistics 33, 307-335.

\bibitem{smith} Smith, R.L. (1990). Max-stable processes and spatial extremes.
Preprint, Univ. North Carolina, USA.

\bibitem{smith+weissman} Smith, R.L., Weissman, I. (1996). Characterization and estimation of the
multivariate extremal index. Manuscript, UNC.

\bibitem{sib} {Sibuya, M.} (1960). Bivariate extreme statistics. \emph{Ann.\hs Inst.\hs Statist.\hs Math.} {\bf 11,} 195-210.

\bibitem{tiago} Tiago de Oliveira, J. (1962/63). Structure theory of bivariate
extremes, extensions. Est. Mat., Estat. e Econ. 7, 165-195.

\bibitem{wolff} Wolff, E. F. (1980). N-dimensional measures of dependence.
Stochastica 4 (3), 175-188.

\end{thebibliography}
\end{document}